\documentclass[conference,a4paper,dvipdfm]{IEEEtran}
\usepackage{cite}

\ifCLASSINFOpdf
\else
\fi
\usepackage[cmex10]{amsmath}
\usepackage{array}

\usepackage{mdwmath}
\usepackage{mdwtab}
\usepackage{fixltx2e}

\usepackage{stfloats}

\usepackage{url}

\usepackage[american,fulldiode]{circuitikz} 
\usepackage{graphicx} 
\usepackage{multirow}
\usepackage{placeins}
\usepackage{amsfonts}
\usepackage{bbm}

\usepackage{slashbox}
\usepackage{amssymb}

\usepackage{lipsum,amsmath,multicol}

\newcommand{\matpower}{M\textsc{ATPOWER}}
\usepackage{multirow}

\let\oldbibliography\thebibliography
\renewcommand{\thebibliography}[1]{%
  \oldbibliography{#1}%
  \setlength{\itemsep}{0pt}%
}

\begin{document}
%
\title{Moment-Based Relaxation of the Optimal Power Flow Problem}

\author{\IEEEauthorblockN{Daniel K. Molzahn, Ian A. Hiskens}
\IEEEauthorblockA{University of Michigan\\
Ann Arbor, MI 48109, USA \\
molzahn@umich.edu, hiskens@umich.edu}
}%


%


\maketitle

\begin{abstract}
The optimal power flow (OPF) problem minimizes power system operating cost subject to both engineering and network constraints. With the potential to find global solutions, significant research interest has focused on convex relaxations of the non-convex AC OPF problem. This paper investigates ``moment-based'' relaxations of the OPF problem developed from the theory of polynomial optimization problems. At the cost of increased computational requirements, moment-based relaxations are generally tighter than the semidefinite relaxation employed in previous research, thus resulting in global solutions for a broader class of OPF problems. Exploration of the feasible space for test systems illustrates the effectiveness of the moment-based relaxation.
\end{abstract}


\vspace{3pt}
\textbf{\textit{Index Terms--Optimal power flow, Global optimization, Moment relaxation, Semidefinite programming}}

%
\IEEEpeerreviewmaketitle

\section{Introduction}
\label{l:intro}

The optimal power flow (OPF) problem determines an optimal operating point for an electric power system in terms of a specified objective function, subject to both network equality constraints (i.e., the power flow equations, which model the relationship between voltages and power injections) and engineering limits (e.g., inequality constraints on voltage magnitudes, active and reactive power generations, and line flows). Generation cost per unit time is typically chosen as the objective function.

The OPF problem is generally non-convex due to the non-linear power flow equations~\cite{bernie_opfconvexity} and may have local solutions~\cite{bukhsh_tps}. Non-convexity of the OPF problem has made solution techniques an ongoing research topic. Many OPF solution techniques have been proposed, including successive quadratic programs, Lagrangian relaxation, genetic algorithms, particle swarm optimization, and interior point methods~\cite{frank2012opf}.

Recently, significant attention has focused on a semidefinite relaxation of the OPF problem~\cite{lavaei_tps}. Using a rank relaxation, the OPF problem is reformulated as a convex semidefinite program. If the relaxed problem satisfies a rank condition (i.e., the relaxation is said to be ``exact'' or ``tight''), the global solution to the original OPF problem can be determined in polynomial time. Prior OPF solution methods do not guarantee finding a global solution in polynomial time; semidefinite programming approaches thus have a substantial advantage over traditional solution techniques. However, the rank condition is not satisfied for all practical OPF problems~\cite{bukhsh_tps,allerton2011,hicss2014}. This paper presents  alternative ``moment-based'' relaxations that globally solve a broader class of OPF problems.


Currently, there is substantial research interest in determining sufficient conditions for which the semidefinite relaxation of~\cite{lavaei_tps} is exact. Existing sufficient conditions include requirements on power injection and voltage magnitude limits and either radial networks (typical of distribution system models) or appropriate placement of controllable phase shifting transformers. (See~\cite{low_irep2013} and the references therein for detailed descriptions of these conditions.) 

Extending this literature to mesh networks without phase-shifting transformers, \cite{madani_asilomar2013} investigates the feasible space of active power injections for weakly-cyclic networks (i.e., networks where no line belongs to more than one cycle). Alternative representations of the line-flow constraints (apparent power, active power, voltage difference, and angle difference limits), while similar in the original OPF problem, can lead to significantly different results for the semidefinite relaxation. Assuming no lower limits on active and reactive power injections at each bus and with line-flow limits represented as voltage differences (i.e., for voltage phasors $V_i$ and $V_j$, constrain the flow between buses $i$ and $j$ as $\left|V_i - V_j \right| \leq \Delta V_{ij}^{\max}$), the semidefinite relaxation of~\cite{lavaei_tps} is proven exact for weakly-cyclic networks with cycles of size three.

While the sufficient conditions developed thus far are promising, they only apply to a limited subset of OPF problems. For more general cases,~\cite{bitar_allerton2013} proposes a method for finding a globally optimal solution that is ``hidden'' in a higher-rank subspace of solutions to the semidefinite relaxation. That is, the solution to the semidefinite relaxation obtained numerically does not satisfy the rank condition but a rank one solution exists with the same globally optimal objective value. For these cases, the semidefinite relaxation is exact in the sense that it yields the globally optimal objective value rather than a strict lower bound, but this fact is not evident as the semidefinite relaxation solution does not directly provide globally optimal decision variables (i.e., the optimal voltage phasors). With a heuristic for finding such ``hidden'' solutions, \cite{bitar_allerton2013} broadens the applicability of the semidefinite relaxation.

However, there exist practical problems for which the semidefinite relaxation is not exact (i.e., the semidefinite relaxation solution has optimal objective value strictly less than the global minimum)~\cite{bitar_allerton2013,hicss2014}. For such cases, \cite{madani_asilomar2013} and~\cite{bitar_allerton2013} propose heuristics for obtaining an (only-guaranteed-locally) optimal solution from the semidefinite relaxation. Heuristics are promising for finding local solutions, and the optimal objective value of the semidefinite relaxation provides a metric for the potential suboptimality of these solutions. For example, the heuristic in~\cite{bitar_allerton2013} finds a solution to a modified form of the IEEE~14-bus system that is within 0.13\% of global optimality. This compares favorably to a solution from the interior-point solver in \matpower~\cite{matpower}, which is only within 4.83\% of global optimality.

While deserving of further study, heuristics eliminate the global optimality guarantee that is one of the main advantages of the semidefinite relaxation. This paper therefore proposes an alternative moment-based convex relaxation that, when exact, yields the global optimum. Using theory developed for polynomial optimization problems~\cite{lasserre_book}, moment-based relaxations have the potential to globally solve a broad class of OPF problems, including many problems for which the semidefinite relaxation of~\cite{lavaei_tps} is not exact. The moment-based relaxation exploits the fact that the OPF problem is composed of polynomials in the voltage phasor components and is therefore a polynomial optimization problem.

The ability to globally solve a broader class of OPF problems has a computational cost. Whereas the semidefinite relaxation of~\cite{lavaei_tps} optimizes matrices composed of all degree-two combinations of the voltage phasor components, the moment-based relaxations optimize matrices composed of higher-degree combinations. In particular, for an $n$-bus system, the order-$\gamma$ moment-based relaxation is solved using a semidefinite program which has a positive semidefinite constraint on a $k\times k$ matrix, where $k = \left(2n+\gamma\right)! / \left( \left(2n\right)! \gamma!\right)$ (i.e., this matrix is composed of all combinations of voltage components up to order $2\gamma$). Thus, the computational requirements of the moment-based relaxations can be substantially larger than the semidefinite relaxation of~\cite{lavaei_tps}, especially for higher orders of the moment-based relaxation.

Fortunately, experience with small systems suggests that low (often second) order relaxations globally solve a broad class of OPF problems, including problems for which the semidefinite relaxation of~\cite{lavaei_tps} is not exact due to disconnected or otherwise non-convex feasible spaces. Note that no fixed-order relaxation is exact for all OPF problems due to the polynomial-time complexity of semidefinite programs as compared to the NP-hardness of some OPF problems~\cite{lavaei_tps}. Indeed, as discussed in Section~\ref{l:results}, some of the NP-hard problems in~\cite{lavaei_tps} provide examples where low-order moment-based relaxations are not exact. Also note that the moment-based relaxation is currently only computationally tractable for small OPF problems; future work includes exploiting sparsity to extend the moment-based relaxations to larger systems.


After introducing the OPF problem formulation in Section~\ref{l:opf_formulation} and describing the moment-based relaxations in Section~\ref{l:msdp_overview}, we explore the feasible space of the second-order moment-based relaxation for a two-bus system in Section~\ref{l:twobusillustration}. For some choices of parameters for this system, the semidefinite relaxation of~\cite{lavaei_tps} is not exact. Since, conversely, the moment-based relaxation is exact for this problem, a comparison of the feasible spaces of the relaxations illustrates the effectiveness of the proposed approach. Section~\ref{l:results} then presents results from the application of the moment-based relaxation to other small OPF problems for which the semidefinite relaxation of~\cite{lavaei_tps} is not exact. Section~\ref{l:conclusion} concludes the paper and discusses future research directions.

\section{OPF Problem Formulation}
\label{l:opf_formulation}

We first present the OPF problem as it is classically formulated. This formulation is in terms of rectangular voltage coordinates, active and reactive power generation, and apparent-power line-flow limits. Consider an $n$-bus power system, where $\mathcal{N} = \left\lbrace 1, 2, \ldots, n \right\rbrace$ is the set of all buses, $\mathcal{G}$ is the set of generator buses, and $\mathcal{L}$ is the set of all lines. $P_{Dk} + j Q_{Dk}$ represents the active and reactive load demand at each bus $k \in \mathcal{N}$. $V_k = V_{dk} + j V_{qk}$ represents the voltage phasors in rectangular coordinates at each bus $k \in \mathcal{N}$. Superscripts ``max'' and ``min'' denote specified upper and lower limits. Buses without generators have maximum and minimum generation set to zero (i.e., $P_{Gk}^{\max} = P_{Gk}^{\min} = Q_{Gk}^{\max} = Q_{Gk}^{\min} = 0, \;\; \forall k\in \mathcal{N}\setminus\mathcal{G}$). $\mathbf{Y} = \mathbf{G} + j \mathbf{B}$ denotes the network admittance matrix.


The network physics are described by the power flow equations:

\begin{subequations}
\small
\begin{align}\nonumber
P_{Gk} = & f_{Pk}\left(V_d,V_q\right) = V_{dk} \sum_{i=1}^n \left( \mathbf{G}_{ik} V_{di} - \mathbf{B}_{ik} V_{qi} \right) &  &  \\ 
\label{opf_Pbalance}  & + V_{qk} \sum_{i=1}^n \left( \mathbf{B}_{ik}V_{di} + \mathbf{G}_{ik}V_{qi} \right) + P_{Dk}  \\ \nonumber
Q_{Gk} = & f_{Qk}\left(V_d,V_q \right) = V_{dk} \sum_{i=1}^n \left( -\mathbf{B}_{ik}V_{di} - \mathbf{G}_{ik} V_{qi}\right) \\
\label{opf_Qbalance} & + V_{qk} \sum_{i=1}^n \left( \mathbf{G}_{ik} V_{di} - \mathbf{B}_{ik} V_{qi}\right) + Q_{Dk}
\end{align}
\end{subequations}

Define a convex quadratic cost function for active power generation:

\begin{equation}\label{objfunction}
f_{Ck}\left(V_d,V_q\right) = c_{k2} \left(f_{Pk}\left(V_d,V_q\right)\right)^2 + c_{k1} f_{Pk}\left(V_d,V_q\right) + c_{k0}
\end{equation}

Define a function for squared voltage magnitude:

\begin{equation} \label{opf_Vsq}
\left(V_{k}\right)^2 = f_{Vk}\left(V_d, V_q\right) = V_{dk}^2 + V_{qk}^2
\end{equation}


Squared apparent-power line-flows $\left(S_{lm}\right)^2$ are polynomial functions of the voltage components $V_d$ and $V_q$. We assume a $\pi$-model with series admittance $g_{lm} + j b_{lm}$ and total shunt susceptance $b_{sh,lm}$ for the line from bus~$l$ to bus~$m$. (For inductive lines and capacitive shunt susceptances, $b_{lm}$ is a negative quantity and $b_{sh,lm}$ is a positive quantity.)

\begin{subequations}
\small
\begin{align}
\nonumber  & P_{lm} = f_{Plm}\left(V_d,V_q\right) = b_{lm} \left( V_{dl} V_{qm} - V_{dm} V_{ql} \right) \\
\label{Plm} & \quad + g_{lm} \left(V_{dl}^2 + V_{ql}^2 - V_{ql} V_{qm} - V_{dl} V_{dm} \right) \\
\nonumber & Q_{lm} = f_{Qlm}\left(V_d,V_q\right)  = b_{lm}\left(V_{dl} V_{dm} + V_{ql} V_{qm} - V_{dl}^2 - V_{ql}^2 \right) \\ 
\label{Qlm} & \quad + g_{lm} \left( V_{dl} V_{qm} - V_{dm} V_{ql} \right) - \frac{b_{sh,lm}}{2} \left( V_{dl}^2 + V_{ql}^2\right) \\
\label{Slm} & \left(S_{lm}\right)^2 = f_{Slm}\left(V_d,V_q\right) = \left(f_{Plm}\left(V_d,V_q\right)\right)^2 + \left(f_{Qlm}\left(V_d,V_q\right)\right)^2
\end{align}
\end{subequations}

\begin{figure*}[!b]
\small
\setcounter{equation}{9}
\vspace*{1pt}
\hrule
\vspace*{0pt}
\begin{equation}\label{2busMomentX} \small
x_2 = \left[\begin{array}{cccccccccc} 1 & V_{d1} & V_{d2} & V_{q2} & V_{d1}^2 & V_{d1}V_{d2} & V_{d1}V_{q2} & V_{d2}^2 & V_{d2}V_{q2} & V_{q2}^2\end{array}\right]^\intercal
\end{equation}
\begin{equation}\label{2busMomentMat} \small
\mathbf{M}_2 \left(y \right) = L_y\left(x_2 x_2^\intercal\right) = \left[\begin{array}{c|ccc|cccccc} 
y_{000} & y_{100} & y_{010} & y_{001} & y_{200} & y_{110} & y_{101} & y_{020} & y_{011} & y_{002} \\\hline
y_{100} & y_{200} & y_{110} & y_{101} & y_{300} & y_{210} & y_{201} & y_{120} & y_{111} & y_{102} \\
y_{010} & y_{110} & y_{020} & y_{011} & y_{210} & y_{120} & y_{111} & y_{030} & y_{021} & y_{012} \\
y_{001} & y_{101} & y_{011} & y_{002} & y_{201} & y_{111} & y_{102} & y_{021} & y_{012} & y_{003} \\ \hline
y_{200} & y_{300} & y_{210} & y_{201} & y_{400} & y_{310} & y_{301} & y_{220} & y_{211} & y_{202} \\ 
y_{110} & y_{210} & y_{120} & y_{111} & y_{310} & y_{220} & y_{211} & y_{130} & y_{121} & y_{112} \\
y_{101} & y_{201} & y_{111} & y_{102} & y_{301} & y_{211} & y_{202} & y_{121} & y_{112} & y_{103} \\
y_{020} & y_{120} & y_{030} & y_{021} & y_{220} & y_{130} & y_{121} & y_{040} & y_{031} & y_{022} \\
y_{011} & y_{111} & y_{021} & y_{012} & y_{211} & y_{121} & y_{112} & y_{031} & y_{022} & y_{013} \\ 
y_{002} & y_{102} & y_{012} & y_{003} & y_{202} & y_{112} & y_{103} & y_{022} & y_{013} & y_{004} 
\end{array}\right]
\end{equation}
\setcounter{equation}{4}
\end{figure*}

The classical OPF problem is then





\begin{subequations}
\label{opf}
\small
\begin{align}
\label{opf_obj} & \min_{V_d,V_q} \sum_{k \in \mathcal{G}} f_{Ck}\left(V_d,V_q\right) \qquad \mathrm{subject\; to} \hspace{-20pt} & \\
\label{opf_P} &  P_{Gk}^{\mathrm{min}} \leq f_{Pk}\left(V_d,V_q \right) \leq P_{Gk}^{\mathrm{max}} & \forall k \in \mathcal{N} \\
\label{opf_Q} &  Q_{Gk}^{\mathrm{min}} \leq f_{Qk}\left(V_d,V_q \right) \leq Q_{Gk}^{\mathrm{max}} &  \forall k \in \mathcal{N} \\
\label{opf_V} &  \left(V_{k}^{\mathrm{min}}\right)^2 \leq f_{Vk}\left(V_d, V_q\right) \leq \left(\vphantom{V_{k}^{\mathrm{min}}} V_{k}^{\mathrm{max}}\right)^2 &  \forall k \in \mathcal{N}  \\
\label{opf_Slm} &  f_{Slm}\left(V_d,V_q\right) \leq \left(S_{lm}^{\mathrm{max}}\right)^2 &  \forall \left(l,m\right) \in \mathcal{L} \\
\label{opf_Vref} & V_{q1} = 0
\end{align}
\end{subequations}

\noindent Note that this formulation limits the apparent-power flow measured at each end of a given line, recognizing that line losses can cause these quantities to differ. Constraint~\eqref{opf_Vref} sets the reference bus angle to zero.

\section{Moment-Based Relaxation Overview}
\label{l:msdp_overview}

The OPF problem~\eqref{opf} is comprised of polynomial functions of the voltage components $V_d$ and $V_q$ and can therefore be solved using moment-based relaxations~\cite{lasserre_book}. We next present the moment-based relaxation for the OPF problem~\eqref{opf}. More detailed descriptions of moment-based relaxations are available in~\cite{lasserre_book}

\begin{figure*}
\small
\setcounter{equation}{12}
\begin{align}\label{2busLocalizingMat} \small \mathbf{M}_{1}\left\lbrace\left(f_{V2} - 0.9\right) y \right\rbrace = \begin{bmatrix}
y_{020} + y_{002} - 0.9y_{000} & y_{120} + y_{102} - 0.9y_{100} & y_{030} + y_{012} - 0.9y_{010} & y_{021} + y_{003} - 0.9y_{001}  \\
y_{120} + y_{102} - 0.9y_{100} & y_{220} + y_{202} - 0.9y_{200} & y_{130} + y_{112} - 0.9y_{110} & y_{121} + y_{103} - 0.9y_{101} \\
y_{030} + y_{012} - 0.9y_{010} & y_{130} + y_{112} - 0.9y_{110} & y_{040} + y_{022} - 0.9y_{020} & y_{031} + y_{013} - 0.9y_{011} \\
y_{021} + y_{003} - 0.9y_{001} & y_{121} + y_{103} - 0.9y_{101} & y_{031} + y_{013} - 0.9y_{011} & y_{022} + y_{004} - 0.9y_{002}  \\
\end{bmatrix}
\end{align}
\setcounter{equation}{5}
\vspace*{1pt}
\hrule
\vspace*{0pt}
\end{figure*}

Polynomial optimization problems, such as the OPF problem, are a special case of a class of problems known as ``generalized moment problems''~\cite{lasserre_book}. Global solutions to generalized moment problems can be approximated using moment-based relaxations that are formulated as semidefinite programs. For polynomial optimization problems that satisfy a technical condition on the compactness of at least one constraint polynomial, the approximation converges to the global solution(s) as the relaxation order increases~\cite{lasserre_book}. (This technical condition can always be satisfied by adding large bounds on all variables and is therefore not restrictive for OPF problems.) Note that while moment-based relaxations can find all global solutions to polynomial optimization problems, we focus on problems with a single global optimum.

Formulating the moment-based relaxation requires several definitions. Define the vector $\hat{x} = \begin{bmatrix} V_{d1} & V_{d2} & \ldots & V_{qn} \end{bmatrix}^\intercal$, which contains all first-order monomials. Given a vector $\alpha = \begin{bmatrix} \alpha_1 & \alpha_2 & \ldots & \alpha_{2n} \end{bmatrix}^\intercal$ with $\alpha \in \mathbb{N}^{2n}$ representing monomial exponents, the expression $\hat{x}^\alpha = V_{d1}^{\alpha_1}V_{d2}^{\alpha_2}\cdots V_{qn}^{\alpha_{2n}}$ defines the monomial associated with $\hat{x}$ and $\alpha$. A polynomial $g\left(\hat{x}\right)$ is then

\begin{equation}\label{gpoly}
g\left(\hat{x}\right) = \sum_{\alpha \in \mathbb{N}^{2n}} g_{\alpha} \hat{x}^{\alpha}
\end{equation}

\noindent where $g_{\alpha}$ is the scalar coefficient corresponding to $\hat{x}^{\alpha}$.

Next define a linear functional $L_y\left\lbrace g\right\rbrace$. 

\begin{equation}\label{L}
L_y\left\lbrace g\right\rbrace = \sum_{\alpha \in \mathbb{N}^{2n}} g_{\alpha} y_{\alpha}
\end{equation}

\noindent This functional replaces the monomials $\hat{x}^\alpha$ in a polynomial function $g\left(\hat{x}\right)$ with scalar variables $y_{\alpha}$. If the argument to the functional $L_y\left\lbrace g \right\rbrace$ is a matrix, the functional is applied to each element of the matrix.

Consider, for example, the vector $\hat{x} = \begin{bmatrix}V_{d1} & V_{d2} & V_{q2} \end{bmatrix}^\intercal$ corresponding to the voltage components of a two-bus system. (For notational convenience, the angle reference constraint $V_{q1}=0$ is used to eliminate $V_{q1}$.) Consider the polynomial $g\left(\hat{x}\right) = -1 + f_{V2}\left(V_d,V_q \right) = -1 + V_{d2}^2 + V_{q2}^2$. (The constraint $g\left(\hat{x}\right) = 0$ forces the squared voltage magnitude at bus~2 to equal 1 per unit.) Then $L_y\left\lbrace g\right\rbrace = -y_{000} + y_{020} + y_{002}$. Thus, $L\left\lbrace g \right\rbrace$ converts a polynomial $g\left(\hat{x}\right)$ to a linear function of $y$. 

The order-$\gamma$ moment-based relaxation forms a vector $x_\gamma$ composed of all monomials of the voltage components up to order $\gamma$:

\begin{align} \nonumber
x_\gamma = & \left[ \begin{array}{ccccccc} 1 & V_{d1} & \ldots & V_{qn} & V_{d1}^2 & V_{d1}V_{d2} & \ldots \end{array} \right. \\ \label{x_d}
& \qquad \left.\begin{array}{cccccc} \ldots & V_{qn}^2 & V_{d1}^3 & V_{d1}^2 V_{d2} & \ldots & V_{qn}^\gamma \end{array}\right]^\intercal
\end{align}

We now define \emph{moment} and \emph{localizing} matrices. The moment matrix $\mathbf{M}_\gamma\left(y\right)$ has entries $y_\alpha$ corresponding to all monomials $x^\alpha$ up to order $2\gamma$ and is symmetric,

\begin{equation}\label{momentmat}
\mathbf{M}_\gamma \left( y \right) = L_y\left( x_\gamma x_\gamma^\intercal\right).
\end{equation}

\noindent Consider, for instance, a two-bus example system with $x_2$ given in~\eqref{2busMomentX}. For $\gamma=2$, this system has the moment matrix shown in~\eqref{2busMomentMat}. This matrix has entries $y_\alpha$ corresponding to all monomials $x^\alpha$ up to degree four.

Note that several terms are repeated in the moment matrix beyond those expected for a generic symmetric matrix. In $\mathbf{M}_2\left(y\right)$, for instance, the terms corresponding to second-order monomials (e.g., $y_{200}$) appear in both the second diagonal block of $\mathbf{M}_2\left(y\right)$ and the first row and column. There are also repetitions in the off-diagonal block, whose entries correspond to third-order monomials (e.g., $y_{210}$) and in the third diagonal block of $\mathbf{M}_2$, whose entries correspond to fourth-order monomials (e.g., $y_{220}$). These repetitions require equality constraints in the semidefinite program implementation.

Symmetric localizing matrices are defined for each constraint of~\eqref{opf}. The localizing matrices consist of linear combinations of the moment matrix entries $y$. Each polynomial constraint of the form $f - a \geq 0$ in~\eqref{opf} (e.g., $f_{V2} - V_2^{\min} \geq 0$) corresponds to the localizing matrix

\setcounter{equation}{11}
\begin{equation}\label{localizing} \small
\mathbf{M}_{\gamma-\beta}\left\lbrace \left(f - a\right) y \right\rbrace = L_y\left\lbrace\left(f - a\right) x_{\gamma-\beta} x_{\gamma-\beta}^\intercal \right\rbrace
\end{equation}
\setcounter{equation}{13}
\vspace{-10pt}

\noindent where the polynomial $f$ has degree $2\beta$ or $2\beta-1$. The bus~2 lower voltage limit with $V_2^{\min} = 0.9$ per unit for the two-bus example system, for example, has the corresponding localizing matrix in~\eqref{2busLocalizingMat}.

We can now form the order-$\gamma$ moment-based relaxation of the OPF problem.

\vspace{-5pt}
\begin{subequations}\small
\label{msdp_opf}
\begin{align}
\label{msdp_obj}& \min_{y} L_y\left\lbrace \sum_{k \in \mathcal{G}} f_{Ck} \right\rbrace \qquad \mathrm{subject\; to} \hspace{-150pt} &  \\
\label{msdp_Pmin} &  \mathbf{M}_{\gamma-1}\left\lbrace \left(f_{Pk} - P_k^{\min}\right) y \right\rbrace \succeq 0 & \forall k\in\mathcal{N}\\
\label{msdp_Pmax} &  \mathbf{M}_{\gamma-1}\left\lbrace \left(P_k^{\max} - f_{Pk} \vphantom{P_k^{\min}}\right) y \right\rbrace \succeq 0 & \forall k\in\mathcal{N}\\
\label{msdp_Qmin} &  \mathbf{M}_{\gamma-1}\left\lbrace \left(f_{Qk} - Q_k^{\min}\right) y \right\rbrace \succeq 0 & \forall k\in\mathcal{N}\\
\label{msdp_Qmax} &  \mathbf{M}_{\gamma-1}\left\lbrace \left(Q_k^{\max} - f_{Qk}  \vphantom{P_k^{\min}}\right) y \right\rbrace \succeq 0 & \forall k\in\mathcal{N}\\
\label{msdp_Vmin} &  \mathbf{M}_{\gamma-1}\left\lbrace \left(f_{Vk} - V_k^{\min}\right) y \right\rbrace \succeq 0 & \forall k\in\mathcal{N}\\
\label{msdp_Vmax} &  \mathbf{M}_{\gamma-1}\left\lbrace \left(V_k^{\max} - f_{Vk}  \vphantom{P_k^{\min}}\right) y \right\rbrace \succeq 0 & \forall k\in\mathcal{N} \\
\label{msdp_Smax} & \mathbf{M}_{\gamma-2}\left\lbrace \left(S_{lm}^{\max} - f_{Slm} \vphantom{P_k^{\min}}\right) y \right\rbrace \succeq 0 & \forall \left(l,m\right)\in\mathcal{L} \\
\label{msdp_Msdp} & \mathbf{M}_\gamma \left(y\right) \succeq 0 & \\
\label{msdp_y0} & y_{00\ldots 0} = 1 & \\
\label{msdp_Vref} & y_{\star\eta\star\star\ldots\star} = 0 & \forall \eta \geq 1
\end{align}
\end{subequations}

\noindent where $\succeq 0$ indicates that the corresponding matrix is positive semidefinite and $\star$ represents any integer in $\left[ 0,\; \gamma \right]$. The moment-based relaxation is thus a semidefinite program. (The dual form of the moment-based relaxation is a sum-of-squares program; see~\cite{lasserre_book} for further details on the dual formulation.) Note that the constraint~\eqref{msdp_y0} enforces the fact that $x^{0} = 1$. The constraint~\eqref{msdp_Vref} corresponds to~\eqref{opf_Vref}; the angle reference can alternatively be used to eliminate all terms corresponding to $V_{q1}$ to reduce the size of the semidefinite program.

Equality constraints are modeled as two inequality constraints. Since, for a symmetric matrix $\mathbf{A}$, the constraints $\mathbf{A} \succeq 0$ and $-\mathbf{A} \succeq 0$ imply $\mathbf{A} = 0$, all entries of a localizing matrix corresponding to an equality constrained polynomial (e.g., the power flow constraints at load buses) are zero.

The order-$\gamma$ moment-based relaxation yields a single global solution if $\mathrm{rank}\left(\mathbf{M}_{\gamma}\left(y\right)\right) = 1$. A solution $x^\ast$ to the OPF problem~\eqref{opf} can then be directly determined from the elements of $y$ corresponding to the linear monomials (e.g., $V_{d1}^\ast = y_{100\cdots0}$). If $\mathrm{rank}\left(\mathbf{M}_{\gamma}\left(y\right) \right) > 1$, then there are either multiple global solutions requiring the solution extraction procedure in~\cite{lasserre_book} or the order-$\gamma$ moment-based relaxation is not exact and yields only a lower bound on the optimal objective value. If the order-$\gamma$ moment-based relaxation is not exact, the order-$\left(\gamma + 1\right)$ moment relaxation will improve the lower bound and may give a global solution.


Note that the order $\gamma$ of the moment-based relaxation must be greater than or equal to half of the degree of any polynomial in the OPF problem~\eqref{opf}. All polynomials can then be written as linear functions of the entries of $\mathbf{M}_\gamma$. For instance, the OPF problem with a linear cost function and without apparent-power line-flow limits requires $\gamma \geq 1$. Although direct implementation of~\eqref{opf} requires $\gamma \geq 2$ due to the fourth-order polynomials in the cost function and apparent-power line-flow limits, these limits can be rewritten as second-order polynomials using a Schur complement~\cite{lavaei_tps}. The OPF problem reformulated in this way only requires $\gamma \geq 1$, but experience suggests that direct implementation of~\eqref{opf} has numerical advantages.

It is interesting to compare the first-order moment-based relaxation and the semidefinite relaxation of~\cite{lavaei_tps}. The first-order moment-based relaxation has a moment matrix with elements corresponding to all monomials up to second-order. The localizing matrices are the constraints multiplied by the scalar~1; the localizing matrix constraints for the first-order moment-based relaxation simply enforce the linear scalar equations $L_y\left\lbrace f - a\right\rbrace = \sum_{\alpha \in \mathbb{N}^{2n}} \left( f_{\alpha} y_{\alpha} \right) - y_{0} a \geq 0$ for all constraints in~\eqref{opf}. Thus, the first-order moment-based relaxation is closely related to the semidefinite relaxation of~\cite{lavaei_tps} which has a feasible space defined by a positive semidefinite matrix constraint and linear constraints on the matrix elements. The slight difference between the formulations is that the semidefinite relaxation of~\cite{lavaei_tps} uses a matrix corresponding to only the second-order monomials (e.g., the second diagonal block in~\eqref{2busMomentMat}), whereas the first-order moment-based relaxation additionally has elements corresponding to constant and linear polynomials in its moment matrix $\mathbf{M}_1\left(y\right)$. This may yield different results for the first-order moment-based relaxation and the semidefinite relaxation of~\cite{lavaei_tps} when the relaxations are not exact.



\section{Application to a Two-Bus Example System}
\label{l:twobusillustration}

With three degrees of freedom $V_{d1}$, $V_{d2}$, and $V_{q2}$ (the angle reference constraint~\eqref{opf_Vref} forces $V_{q1} = 0$), the two-bus example system in~\cite{bukhsh2011} allows for visualizing the entire feasible space of the OPF problem~\eqref{opf}. For some choices of parameters, the OPF problem~\eqref{opf} for this system has a disconnected feasible space, and the semidefinite relaxation of~\cite{lavaei_tps} is not exact. (Note that this system does not satisfy the sufficient conditions for exactness of the semidefinite relaxation described in~\cite{low_irep2013}.) This section illustrates the feasible space for the second-order moment-based relaxation, which finds a global solution for a much larger range of parameters for this problem.

\begin{figure}[tbh]
\centering
\includegraphics[totalheight=0.07\textheight]{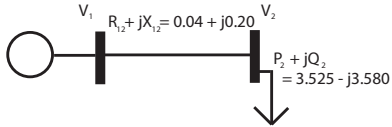}
\caption{Two-Bus System from ~\cite{bukhsh2011}}
\label{f:twobussystem}
\end{figure}

Fig.~\ref{f:twobussystem} gives the system's one-line diagram assuming a 100~MVA base power.  The generator at bus~1 has no limits on active or reactive outputs and there is no line-flow limit. Bus~1 voltage magnitude is in the range $\left[0.95, 1.05\right]$ per unit, while bus~2 voltage magnitude is greater than 0.95 per unit and less than the parameter $V_2^{\max}$.

\begin{figure}[t]
\centering
\includegraphics[totalheight=0.20\textheight]{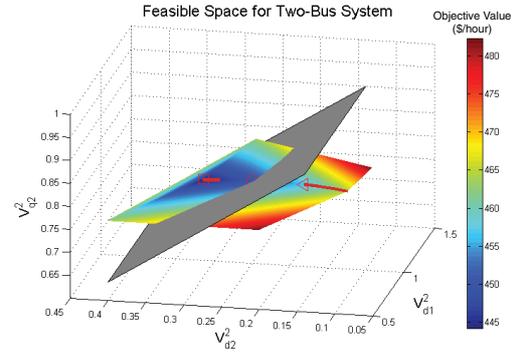}
\caption{Feasible Space of the Semidefinite Relaxation of~\cite{lavaei_tps} for the Two-Bus System Showing the Constraint $V_{2}^{\max} = 1.02$ per unit}
\label{f:case2local_space_cut}
\end{figure}

Fig.~\ref{f:case2local_space_cut} shows the feasible space for the semidefinite relaxation of~\cite{lavaei_tps}. The colored conic shape is the projection of the feasible space of the semidefinite relaxation into the space of squared voltage components ($V_{d1}^2$, $V_{d2}^2$, and $V_{q2}^2$), with the colors based on a \$1/MWh cost of active power generation at bus~1. The red line forms the (disconnected) feasible space for the OPF problem~\eqref{opf}. With $V_{2}^{\max} = 1.05$ per unit, both the semidefinite relaxation of~\cite{lavaei_tps} and the OPF problem~\eqref{opf} have global minimum at the red square in Fig.~\ref{f:case2local_space_cut}, and the semidefinite relaxation is exact. The more stringent limit of $V_{2}^{\max} = 1.02$ per unit is shown by the gray plane; this constraint eliminates the feasible space to the left of the plane. The solution to the semidefinite relaxation of~\cite{lavaei_tps} is at the red circle on the gray plane, while the global solution to the classical OPF problem is at the red triangle. Thus, the semidefinite relaxation of~\cite{lavaei_tps} is not exact. 

\begin{figure}[t]
\centering
\includegraphics[totalheight=0.24\textheight]{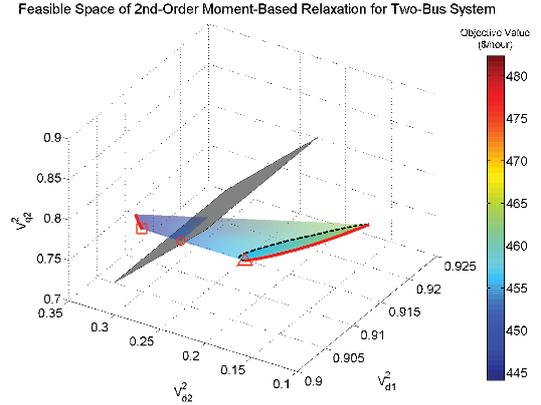}
\caption{Feasible Space of the Second-Order Moment-Based Relaxation for the Two-Bus System Showing the Constraint $V_{2}^{\max} = 1.02$ per unit}
\label{f:case2local_space_msdp_cut}
\end{figure}

Conversely, the second-order moment-based relaxation is exact for both $V_2^{\max} = 1.05$ per unit and $V_2^{\max} = 1.02$ per~unit. Fig.~\ref{f:case2local_space_msdp_cut} shows a projection of the feasible space for this problem. The gray plane again corresponds to $V_{2}^{\max} = 1.02$ per unit in the projected space. The feasible space for the second-order moment-based relaxation is planar with boundaries equal to the feasible space of the OPF problem~\eqref{opf}, which consists solely of the two red line segments on the left and far right of Fig.~\ref{f:case2local_space_msdp_cut}. (Both the colors showing the generation cost and the feasible space values are recovered from the entries of the moment matrix corresponding to the squared terms in the second diagonal block of~\eqref{2busMomentMat}.) With $V_{2}^{\max} = 1.05$ per unit, the second-order moment-based relaxation finds the global solution at the red square in Fig.~\ref{f:case2local_space_msdp_cut}.

In the projection shown in Fig.~\ref{f:case2local_space_msdp_cut}, it appears that imposing the limit $V_2^{\max} = 1.02$ per unit will result in the second-order moment-based relaxation finding the point at the red circle (which does not satisfy the condition $\mathrm{rank}\left(\mathbf{M}_2\left( y \right) \right) = 1$) rather than the global optimum to the OPF problem~\eqref{opf}, which is at the red triangle. (Points on the plane between the red line segments, such as the red circle, are \emph{not} feasible for the OPF problem, but are in the feasible space of the second-order moment-based relaxation with $V_2^{\max} = 1.05$ per unit.)

However, the second-order moment-based relaxation finds the global solution at the red triangle and is therefore exact. While the red circle is in the feasible space when $V_2^{\max} = 1.05$ per unit, this and nearby points are eliminated from the feasible space when $V_2^{\max} = 1.02$ per unit by the upper voltage magnitude constraint~\eqref{msdp_Vmax}. That is, the localizing matrix $\mathbf{M}_{1}\left\lbrace \left(V_2^{\max} - f_{V2} \vphantom{P_k^{\min}}\right) y \right\rbrace$ is positive semidefinite for these points when $V_2^{\max} = 1.05$ per unit but not when $V_2^{\max} = 1.02$ per unit. The feasible space with $V_2^{\max}=1.02$ per unit is the planar region between the black dashed line and the red line on the right. (Note that there is a small range $V_{2}^{\max} \in \left[1.0336, 1.0338 \right]$ per unit where the second-order relaxation does not yield a global optimum; a third-order relaxation finds the global solution with $V_{2}^{\max}$ in this range.)

With the need to have consistent higher-order terms in $y$ that yield positive semidefinite moment and localizing matrices, moment-based relaxations with $\gamma > 1$ are tighter than the semidefinite relaxation of~\cite{lavaei_tps}. The improved tightness has a computational cost: the largest matrix in the semidefinite relaxation of~\cite{lavaei_tps} is $3\times 3$ in contrast to $10 \times 10$ for the second-order moment-based relaxation.


\section{Results for Previously Published Examples}
\label{l:results}

Section~\ref{l:twobusillustration} shows how a moment-based relaxation globally solves a problem for which the semidefinite relaxation of~\cite{lavaei_tps} is not exact. The moment-based relaxation is next applied to other small problems for which the semidefinite relaxation of~\cite{lavaei_tps} is not exact. These problems were solved using YALMIP's moment solver~\cite{yalmip} and SeDuMi~\cite{sedumi}.

Table~\ref{t:results_small} lists small problems for which the semidefinite relaxation of~\cite{lavaei_tps} is not exact for certain parameters. The number of buses is appended to the case names. The table shows the lowest order $\gamma_{\min}$ needed for a global solution. A second-order moment-based relaxation suffices for a broad class of problems. Third-order relaxations are occasionally needed for small parameter ranges.

\begin{table}[tbh]
\centering
\begin{tabular}{|l|l|c|}
\hline 
\textbf{Case} & \textbf{Parameters} & $\gamma_{\min}$\\ \hline\hline
\emph{LMBD3}~\cite{allerton2011} & 50 MVA line limit & 2 \\\hline
\emph{MLD3}~\cite{hicss2014} & 100 MVA line limit & 2 \\\hline
\emph{BGMT3}~\cite{bukhsh_tps} & & 2 \\\hline
\emph{LH5}~\cite{bernie_opfconvexity} & $P_{D3} = 17.17$ per unit & 2 \\\hline
\multirow{3}{40pt}{\emph{BGMT5}~\cite{bukhsh_tps}} & $Q_2^{\min} \in \left[-50, -27.36 \right]$ MVAR & 2 \\
 & $Q_2^{\min} \in \left[-27.35, -27.04 \right]$ MVAR & 3 \\
 & $Q_2^{\min} \in \left[-27.03, 0 \right]$ MVAR & 2 \\ \hline
\emph{BGMT9}~\cite{bukhsh_tps} & & 2 \\\hline
\emph{MSL10\_ex1}~\cite{madani_asilomar2013} &  & $>2$ \\ \hline
\emph{MSL10\_ex2}~\cite{madani_asilomar2013} &  & 2\\\hline
\end{tabular}
\caption{Moment-Based Relaxation Results}
\label{t:results_small}
\vspace{-1.5em}
\end{table}


We have some specific comments on these examples. The three-bus OPF problems \emph{LMBD3} and \emph{MLD3} have binding apparent-power line-flow limits. The line-flow limit in \emph{MLD3} results in a disconnected feasible space~\cite{hicss2014}. Thus, the second-order moment-based relaxation globally solves at least some problems for which the semidefinite relaxation of~\cite{lavaei_tps} is not exact due to tight line-flow limits. 

The semidefinite relaxation of~\cite{lavaei_tps} is not exact for the three and nine-bus problems \emph{BGMT3} and \emph{BGMT9} due to the presence of local optima. The second-order moment-based relaxation finds the global optimum for these problems. \matpower~\cite{matpower} with the MIPS solver initialized using a ``flat start'' (unity voltage magnitudes with zero voltage angles) finds a local optimum for \emph{BGMT9} with objective value 38\% greater than the globally optimal objective value. The semidefinite relaxation of~\cite{lavaei_tps} yields a lower bound that is 11\% less than the global optimum. Thus, existing techniques do not perform well for this problem.

Similar to the problem in Section~\ref{l:twobusillustration}, a third-order moment-based relaxation is needed to globally solve the five-bus problem \emph{BGMT5} for a narrow parameter range.

\begin{figure}[tbh]
\centering
\includegraphics[totalheight=0.27\textheight]{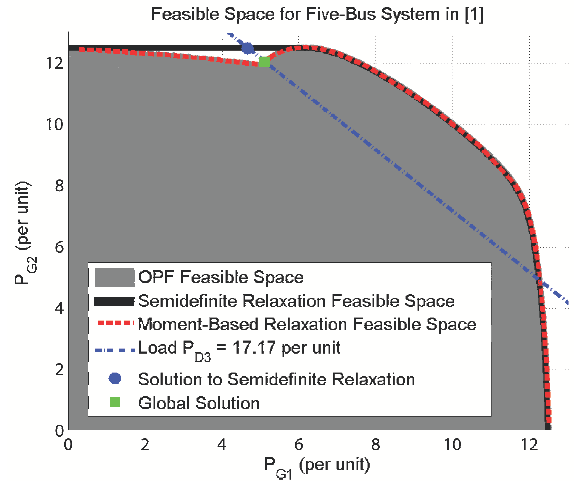}
\caption{Feasible Space for Five-Bus System in~\cite{bernie_opfconvexity}}
\label{f:fivebusspace}
\end{figure}

The gray region in Fig.~\ref{f:fivebusspace} shows the non-convex feasible space of active power injections for the lossless five-bus OPF problem \emph{LH5}. The feasible space of the semidefinite relaxation of~\cite{lavaei_tps} shown by the black curve in Fig.~\ref{f:fivebusspace} is not tight for some parameters (e.g., the dashed blue line representing a load $P_{D3}= 17.17$ per unit with generation $P_{G1}$ more expensive than $P_{G2}$)~\cite{hicss2014}. Conversely, the feasible space for the second-order moment-based relaxation shown by the red dashed curve in Fig.~\ref{f:fivebusspace} is tight for varying load demands $P_{D3}$, which result in other lines parallel to the blue line, and generator costs.



The ten-bus problem \emph{MSL10\_ex1} has at least two global solutions (the solution in~\cite{madani_asilomar2013} and $\begin{bmatrix}P_{G4} \! & \!P_{G5}\! & \!P_{G6}\! & \!P_{G9}\! & \!P_{G10}\end{bmatrix} = \begin{bmatrix}16.99\! & \!32.66\! & \!0\! & \!0\! & \!38.35 \end{bmatrix}$ MW). A moment-based relaxation with $\gamma > 2$ is needed to find the multiple solutions, but only a second-order relaxation is currently computationally tractable.

The second-order moment-based relaxation gives the single global optimum for the related ten-bus problem \emph{MSL10\_ex2}. The heuristic method proposed in~\cite{madani_asilomar2013} finds an only-locally-optimal solution to this problem with an objective value of \$153.97, which is 28.4\% greater than the global minimum of \$119.95 with $\begin{bmatrix}P_{G4} \! & \!P_{G5}\! & \!P_{G6}\! & \!P_{G9}\! & \!P_{G10}\end{bmatrix} = \begin{bmatrix}0\! & \!47.994\! & \!0\! & \!0\! & \!40.006 \end{bmatrix}$ MW.

The moment-based relaxation was also applied to two NP-hard OPF problems from~\cite{lavaei_tps}. The first problem, \emph{LLn\_1} where $n$ represents an arbitrary number of buses, eliminates all active, reactive, and line-flow limits to minimize network losses with voltage magnitudes constrained to 1~per unit. The second problem, \emph{LLn\_2}, minimizes losses for a purely resistive network with zero reactive power injections and voltages constrained to the discrete set $\left\lbrace -1, 1\right\rbrace$. 

Table~\ref{t:nphard} summarizes the application of the moment-based relaxations to these problems. Two network structures were considered: a ring and a complete network (each bus connected to every other bus). In both cases, the global solution was obtained (or the computational capabilities were surpassed without yielding a solution) at the same order of the moment-based relaxation. For \emph{LLn\_1}, lines have 0.1 per unit reactances and either zero or $1\times 10^{-3}$ per unit resistance. For \emph{LLn\_2}, lines have zero reactance and $0.1$ per unit resistance.

\begin{table}[tbh]
\small
\centering
\begin{tabular}{|l|c|c|}
\hline 
\textbf{Case} & \textbf{Parameters} & $\gamma_{\min}$\\ \hline\hline
\emph{LLn\_1}, $n=2,\ldots,9$ & $R = 1\times 10^{-3}$ per unit & 2 \\
\emph{LL2\_1}, & $R = 0$ per unit & $3$ \\
\emph{LL3\_1} & $R = 0$ per unit & $5$ \\
\emph{LL4\_1} & $R = 0$ per unit & $> 4$ \\
\emph{LLn\_1}, $n=5,6$ & $R = 0$ per unit & $> 3$ \\
\emph{LLn\_1}, $n=7,8,9$ & $R = 0$ per unit & $> 2$ \\ \hline
\emph{LLn\_2}, $n=2,\ldots,9$ &  & 2 \\\hline
\end{tabular}
\caption{Moment-Based Relaxation Applied to NP-Hard Problems}
\label{t:nphard}
\vspace{-1.5em}
\end{table}

Table~\ref{t:nphard} shows that low-degree moment-based relaxations have some success for small problems of the forms \emph{LLn\_1} and \emph{LLn\_2}. While a second-order moment-based relaxation solves \emph{LLn\_1} with lossy networks and \emph{LLn\_2}, low-order moment-based relaxations are not well-suited for \emph{LLn\_1} with lossless networks. This is likely due to the fact that \emph{LLn\_1} with lossless networks has multiple global solutions, while \emph{LLn\_1} with lossy networks and \emph{LLn\_2} have unique global solutions.

Since there exist NP-hard OPF problems, the moment-based relaxation for a fixed-order $\gamma$ (or any other polynomial-time relaxation) cannot globally solve all OPF problems. However, the NP-hard problems in~\cite{lavaei_tps} do not represent typical power systems. The results shown in Table~\ref{t:results_small} suggest that the moment-based relaxation may be exact for a broad class of OPF problems that excludes some NP-hard problems.

\section{Conclusions and Future Work}
\label{l:conclusion}

Using theory developed for polynomial optimization problems, this paper has proposed a moment-based relaxation for the OPF problem. With the trade-off of increased computational requirements, the moment-based relaxation globally solves a broader class of problems than previous convex relaxations, such as the semidefinite relaxation of~\cite{lavaei_tps}. After formulating the moment-based relaxation of the OPF problem, this paper has investigated the feasible space of the moment-based relaxation. Global solution of several small problems for which the semidefinite relaxation of~\cite{lavaei_tps} is not exact demonstrates the moment-based relaxation's effectiveness.

There are many avenues for future work on moment-based relaxations, including computational improvements, developing sufficient conditions for exactness of the relaxation, and extension to other power systems problems. The large size of the semidefinite programs used to evaluate the moment-based relaxation currently precludes application to problems with more than ten buses. Future work includes the exploitation of power system sparsity to solve larger OPF problems. Similar to the semidefinite relaxation of~\cite{lavaei_tps}, a matrix completion decomposition~\cite{molzahn_holzer_lesieutre_demarco-large_scale_sdp_opf} is applicable to moment-based relaxations. Further, the results of~\cite{molzahn_holzer_lesieutre_demarco-large_scale_sdp_opf} indicate that only small subnetworks of typical OPF problems cause the semidefinite relaxation of~\cite{lavaei_tps} to not be exact. Therefore, it may be sufficient to selectively apply the moment-based relaxations to small subnetworks of a large OPF problem, thus further improving computational tractability.

Generalizing current efforts to find sufficient conditions for which the semidefinite relaxation of~\cite{lavaei_tps} is exact, future work also includes developing sufficient conditions for which an order-$\gamma$ moment based relaxation is exact.

The moment-based relaxation can also be extended to other power systems problems, including those  problems for which the semidefinite relaxation of~\cite{lavaei_tps} has already shown some success (e.g., state estimation, voltage stability, finding multiple power flow solutions). Further, the ability of moment-based relaxations to include binary variables provides the opportunity to consider problems with discrete constraints such as transmission switching and unit commitment.

\bibliographystyle{IEEEtran}
\bibliography{IEEEabrv,molzahn_hiskens-Moment_OPF_preprint}
%

\end{document}